\documentclass{amsart}

 \newtheorem{thm}{Theorem}[section]
 \newtheorem{cor}[thm]{Corollary}
 \newtheorem{lem}[thm]{Lemma}
 \newtheorem{prop}[thm]{Proposition}
 \theoremstyle{definition}
 \newtheorem{defn}[thm]{Definition}
 \newtheorem{rem}[thm]{Remark}
 \numberwithin{equation}{section}
 \DeclareMathOperator{\rank}{rank}
 \DeclareMathOperator{\Span}{Span}

\begin{document}

\title[Integral Bases for certain TQFT-Modules of the torus]
{Integral Bases for certain TQFT-Modules of the torus}
\author{Khaled Qazaqzeh}

\address{Department of Mathematics, Louisiana State University, Baton Rouge, 70803,
USA} \email{qazaqzeh@math.lsu.edu}
\urladdr{www.math.lsu.edu/\textasciitilde qazaqzeh/}
\date{09/23/2005}


%%% ----------------------------------------------------------------------

\begin{abstract}
We find two bases for the lattices of the $SU(2)$-TQFT-theory
modules of the torus over given rings of integers. We use variant
of the bases defined in ~\cite{GMW04} for the lattices of the
$SO(3)$-TQFT-theory modules of the torus. Moreover, we discuss the
quantization functors ($V_{p},Z_{p}$) for $p = 1$, and $p = 2$. Then
we give concrete bases for the lattices of the modules in the
2-theory. We use the above results to discuss the ideal invariant
defined in ~\cite{FB01}. The ideal can be computed for all the
3-manifolds using the 2-theory, and for all 3-manifolds with torus
boundary using the $SU(2)-$TQFT-theory. In fact, we show that this ideal in
the $SU(2)-$TQFT-theory is contained in the product of the ideals
in the 2-theory and the $SO(3)-$TQFT-theory under a certain change
of coefficients, and it is equal in the case of a torus boundary.
\end{abstract}

%%% ----------------------------------------------------------------------
\maketitle
%%% ----------------------------------------------------------------------

\section*{Introduction}
 We let $p$ denote an odd prime or twice an odd prime
unless mentioned otherwise. Also, we let $\Sigma$ denote a surface
of genus $g$. Gilmer defined an integral TQFT-functor
$\mathcal{S}_{p}$ in ~\cite{G04} based on the integrality results
of the $SO(3)$- and $SU(2)$-invariants in ~\cite{GR97,M}. This is
a functor that associates to a closed surface $\Sigma$, a module
$\mathcal{S}_{p}(\Sigma)$ over a certain cyclotomic ring of
integers $\mathcal{O}_{p}$. Moreover, Gilmer showed that these
modules are free in the case of $p$ is an odd prime. Gilmer and
Masbaum constructed basis for $\mathcal{S}_{p}(\Sigma)$ and gave
an independent proof of freeness in this case. In addition, Gilmer
showed that these modules are projective where $p$ is twice an odd
prime. In this paper, we prove that the modules
$\mathcal{S}_{p}(S^{1}\times S^{1})$ are free by constructing two
explicit bases in the case that $p$ is twice an odd prime. In the
2-theory, we prove also that the modules $\mathcal{S}_{2}(\Sigma)$
are free by constructing an explicit basis for any surface.

Frohman and Kania-Bartoszynska in ~\cite{FB01} defined an ideal
invariant of 3-manifolds with boundary using the
$SU(2)$-TQFT-theory that is hard to compute. In fact, they make
use of another ideal that they defined to give an estimate for
this ideal. However, Gilmer and Masbaum in ~\cite{GM04} computed
an analogous ideal invariant using the $SO(3)-$TQFT-theory for
3-manifolds that are obtained by doing surgery along a knot in the
complement of another knot. The computations depend entirely on
the fact that bases are constructed for the integral lattices of
the $SO(3)$-TQFT-theory modules ~\cite{GMW04,GM04} of the torus.
Also, Gilmer and Masbaum gave a finite set of generators for this
ideal in general. Based on our results in this paper, we compute
this ideal for the above 3-manifolds with torus boundary using the
$SU(2)$-TQFT-theory. Also, we introduce a formula to give an
estimate for the ideal using the $SU(2)$-TQFT-theory in terms of
the ideals using the 2- and $SO(3)$-TQFT-theories. In fact, the
same formula can be used to compute this ideal using the
$SU(2)$-TQFT-theory for the all the above 3-manifolds with torus
boundary.
\par
In $\S1$, we describe the $SO(3)$- and $SU(2)$-TQFT-functors using
the approach of ~\cite{BHMV3} over a variant ring depending on
$p$. We review the integral TQFT-functors in  $\S2$ that Gilmer
defined in \cite{G04}. The first bases for the lattices of the
$SU(2)$-TQFT-modules are given in $\S3$. We review the Frohman
Kania-Bartoszynska ideal in $\S4$, and then we draw some
conclusions based on the results of the previous section regarding
this ideal. The quantization functors for $p=1$ and $p=2$ are
discussed in $\S5$, again following ~\cite{BHMV3}. Also in this
section, we give basis for $\mathcal{S}_{2}(\Sigma)$, and then
draw some conclusions regarding the Frohman and Kania-Bartoszynska
ideal for this theory. We reformulate some results given in
\cite{BHMV3} in $\S6$ to serve our need. Finally, we give another
bases for the lattices of the $SU(2)$-TQFT-modules in $\S7$. The
advantage of this one over the first basis is that it allows us to
prove Theorem (\ref{t:ideal}).

%%% ----------------------------------------------------------------------

\section{The $SO(3)$- and $SU(2)$-TQFTs}
We consider the (2+1)-dimensional TQFT constructed as the main
example of ~\cite[P.\,456]{BHMV3} with some modifications. In
particular, we use the cobordism category $\mathcal{C}$ discussed
in \cite{G04,GQ} where the 3-manifolds have banded links but
surfaces do not have colored points. Hence the objects are
oriented surfaces with extra structure (Lagrangian subspaces of
their first real homology). The cobordisms are equivalence classes
of compact oriented 3-manifolds with extra structure (an integer
weight) with banded links sitting inside of them. Two cobordisms
with the same weight are said to be equivalent if there is an
orientation preserving diffeomorphism that fixes the boundary.

Let
\[ k_{p} =
  \begin{cases}
  \mathbb{Z}[A_{p},\frac{1}{p}],
     & \text{if $p \equiv  \quad -1 \pmod{4}$};
     \\
    {\mathbb{Z}[\alpha_{p},\frac{1}{p}]}, & \text{if $p \equiv 1 $ or $2
    \pmod{4}$}.
  \end{cases}
\]
Here and elsewhere $A_{p}$, $\alpha_{p}$ are $\zeta_{2p}$ and
$\zeta_{4p}$ respectively for $p\geq 3$.

Now, we consider the TQFT-functor ($V_{p},Z_{p}$) from
$\mathcal{C}$  to the category of finitely generated projective
$k_{p}$-modules. The functor ($V_{p},Z_{p}$) is defined as
follows. $V_{p}(\Sigma)$ is a quotient of the $k_{p}$-module
generated by all cobordisms with boundary $\Sigma$, and $Z_{p}(M)$
is the $k_{p}$-linear map from $V_{p}(\Sigma)$ to
$V_{p}(\Sigma^{'})$ (where $\partial M = -\Sigma\coprod
\Sigma^{'}$) induced by gluing representatives of elements of
$V_{p}(\Sigma)$ to $M$ along $\Sigma$ via the identification map
of the first component of the boundary.

If $M$ is a closed cobordism, then $Z_{p}[M]$ is the
multiplication by the scalar $ \langle M \rangle_{p}$ defined in
~\cite[\S.\,2]{BHMV3}. This invariant is normalized in two other
ways. The first normalization of this invariant is
$I_{p}(M)=\mathcal{D}_{p}\langle M_{\flat}\rangle_{p}$. Here and
elsewhere $M_{\flat}$ is the 3-manifold $M$ with a reassigned
weight zero, and $\mathcal{D}_{p}=\langle
S^{3}_{\flat}\rangle_{p}^{-1}$. The second normalization is
$\theta_{p}(M)=\mathcal{D}_{p}^{\beta_{1}(M)}\langle
M_{\flat}\rangle_{p}$, i.e
$\theta_{p}(M)=\mathcal{D}_{p}^{\beta_{1}(M)}I_{p}(M)$.

If $\partial M = \Sigma$ and $M$ is considered as a cobordism from
$\emptyset$ to $\Sigma$, then $Z_{p}(M)(1)\in V_{p}(\Sigma)$ is
denoted by $[M]_{p}$ and called a \textsl{vacuum state} and it is
\textsl{connected} if $M$ is connected. Finally, note that $V_{p}$
is generated over $k_{p}$ by all vacuum states.

 The modules $V_{p}(\Sigma)$ are free modules over $k_{p}$, and carry a nonsingular
 Hermitian bilinear form
 \[
 \langle \ , \ \rangle_{\Sigma} : V_{p}(\Sigma)\times V_{p}(\Sigma)
 \rightarrow k_{p},
 \]
 given by
 \begin{equation}\label{e:form1}
 \langle[M_{1}],[M_{2}]\rangle_{\Sigma} = \langle
 M_{1}\cup_{\Sigma}-M_{2}\rangle_{p}.
 \end{equation}
 Here -$M$ is the cobordism $M$ with the orientation reversed and
 multiplying the integer weight by -1, and leaving the Lagrangian
 subspace on the boundary the same.

  Let $d_{p} = [\frac{p-1}{2}]$, it is known that $d_{p}$ is the dimension of
 $V_{p}(S^{1}\times S^{1})$. One has that $V_{p}(S^{1}\times S^{1})\cong
 k_{p}[z]/I$ where the ideal $I$ is generated by $e_{d_{p}}-e_{d_{p}-1}$
 in the case of $p$ is an odd prime and by $e_{d_{p}}$ in the case of $p$ is
 twice an odd prime (See \cite{BHMV1} for more details).
 Thus indeed, $V_{p}(S^{1}\times S^{1})$ has a
 basis $\{e_{0},\ldots,e_{d_{p}-1}\}$ of rank $d_{p}$.

%%% ----------------------------------------------------------------------

\section{The Integral Cobordism Functor}
Let $\mathcal{C^{'}}$ be the subcategory of $\mathcal{C}$
consisting of the nonempty connected surfaces and connected
cobordisms between them. Let $\mathcal{O}_{p}$ be the ring of
integers of the ring $k_{p}$ defined before. The ring of integers
is given by
 \[ \mathcal{O}_{p} =
  \begin{cases}
  \mathbb{Z}[A_{p}],
     & \text{if $p \equiv  \quad -1 \pmod {4}$};
     \\
    {\mathbb{Z}[\alpha_{p}]}, & \text{if $p \equiv 1 $ or $ 2 \pmod
    {4}$}.
  \end{cases}
\]

Thus the ring of integers of $k_{p}$ is a Dedekind domain.

\begin{defn}
For the surface $\Sigma$, we define $\mathcal{S}_{p}(\Sigma)$ to be the
$\mathcal{O}_{p}$-submodule of $V_{p}(\Sigma)$ generated by all
connected vacuum states.
\end{defn}
If $M :\Sigma \rightarrow \Sigma^{'}$ is a cobordism of
$\mathcal{C}$, then $Z_{p}(M)([N]_{p}) = [M \cup_{\Sigma}-N]_{p}
\in \mathcal{S}_{p}(\Sigma{'})$. Hence we obtain a functor from
$\mathcal{C}^{'}$ to the category of $\mathcal{O}_{p}$-modules.
These modules are projective as they are finitely generated
torsion-free over Dedekind domains \cite[Thm.\,2.5]{G04}. Also, these modules carry
 an $\mathcal{O}_{p}$-Hermitian bilinear form
\[(\ ,\ )_{\Sigma}:\mathcal{S}_{p}(\Sigma)\times \mathcal{S}_{p}(\Sigma) \rightarrow
\mathcal{O}_{p},\] given by
\begin{equation}\label{e:form2}
 ([M_{1}],[M_{2}])_{\Sigma} =
\mathcal{D}_{p}\langle[M_{1}],[M_{2}]\rangle =
\mathcal{D}_{p}\langle M_{1}\cup_{\Sigma}-M_{2}\rangle_{p},
\end{equation}

The value of this form always lies in $\mathcal{O}_{p}$ by the
integrality results for closed 3-manifolds in ~\cite{GR97,M}.

If $R\subseteq \mathcal{S}_{p}(\Sigma)$ is an $\mathcal{O}_{p}-$submodule
define
\[R^{\sharp}=\{v\in V_{p}(\Sigma)| (r,v)_{\Sigma}\in\mathcal{O}_{p}, \forall r\in R\},\]
then we can conclude
\begin{equation}\label{e:dual}
R\subseteq \mathcal{S}_{p}(\Sigma)\subseteq \mathcal{S}_{p}^{\sharp}
(\Sigma)\subseteq R^{\sharp}.
\end{equation}
\begin{defn}
A Hermitian bilinear form on a projective module over a Dedekind domain is
called  $\textit{non-degenerate}$ if the adjoint map is injective,
and $\textit{unimodular}$ if the adjoint map is an isomorphism.
\end{defn}
For our use, if the matrix of the form has a nonzero (unit)
determinant, then the form will be non-degenerate (unimodular)
respectively. Note that the determinant of the form
(\ref{e:form2}) is nonzero from the fact that the form (\ref{e:form1})
is non-degenerate. Hence the form (\ref{e:form2}) is
non-degenerate. In fact, we prove that the form (\ref{e:form2}) is
unimodular for the 2-theory (discussed in $\S5$) for all surfaces and for $S^{1}\times
S^{1}$ in the case of $p$ is twice an odd prime.

A standard basis $\{u_{\sigma}\}$ for $V_{p}(\Sigma)$ is given
(see ~\cite{BHMV3}) in terms of $p$-admissible colorings $\sigma$
of the spine of a handlebody of genus $g$ whose boundary is
$\Sigma$ where the set of colors is $\{0,1,2,\ldots,d_{p}-1\}$,
and the sum of the colors at a 3-vertex even and less than
$2d_{p}$ in the case that $p$ is twice an odd prime.

All of the above elements $u_{\sigma}$ lie in
$\mathcal{S}_{p}(\Sigma)$ when $p$ is twice an odd prime. This
follows from the fact that the quantum integers (denominators of
the Jones-Wenzel idempotents) are units in $\mathcal{O}_{p}$ (see
Corollary \ref{c:unit}). An admissible colored trivalent graph
\cite{BHMV3} is to be interpreted, here and elsewhere, as an
$\mathcal{O}_{p}$-linear combination of links.

We say $a\sim b$ in $\mathcal{O}_{p}$ if $a/b$ is a unit in
$\mathcal{O}_{p}$. The following proposition is an elementary fact
from number theory that gives us a family of units in the ring
$\mathcal{O}_{p}$.
\begin{prop}[\cite{W}]\label{p:unit}
Suppose $n$ has at least two distinct prime factors. Then
$1-\zeta_{n}$ is a unit in $\mathbb{Z}[\zeta_{n}]$.
\end{prop}
We make use of the following lemma in giving the first basis
for $\mathcal{S}_{p}(S^{1}\times S^{1})$ in $\S3$ whose proof will be in $\S6$.
\begin{lem}\label{l:d}
\[ \mathcal{D}_{p} \sim
  \begin{cases}%

    {(1-A_{p}^{2})^{\frac{p-3}{2}}}, & \text{if $p$ is an odd prime};
    \\
    {\sqrt{2}(1+\alpha_{p}^{4})^{\frac{\frac{p}{2}-3}{2}}},&\text{if $p$ is twice an odd
    prime}.
 \end{cases}
\]
\end{lem}
%By proposition (\ref{p:unit}), we conclude that $(1+A^{2})$
%is a unit in $\mathcal{O}_{p}$ in the case of $p$ is an odd prime. We know that
% \[
% \mathcal{D}_{p}^{2}= \frac{-p}{(q-q^{-1})^{2}} \quad \text{where}
% \quad q=A^2.
%\]

%Now the result follows from the fact that
% $p\sim(1-A^{2})^{p-1}$.
%In the other case, the result follows by a similar argument and using the following facts:

%\begin{itemize}
 %   \item  ($1-\alpha^4$) is a unit in $\mathcal{O}_{p}$, where $\alpha^4=q$,
  %  \item  $p/2\sim(1+\alpha^{4})^{\frac{p}{2}-1}$.
%\end{itemize}

\begin{prop}\label{p:pairing}
 The elements $\{u_{\sigma}\}$ are orthogonal with respect to the form
 (\ref{e:form2}). Moreover, we have
 \begin{equation}
 (u_{\sigma},u_{\sigma})_{\Sigma} \sim \mathcal{D}_{p}^{g}.
 \end{equation}

\end{prop}
\begin{proof}
If we use Theorem (4.11) in ~\cite{BHMV3}, and the facts
\[\langle S^{3}\rangle_{p}= \mathcal{D}_{p}^{-1}, \ \sharp
v-\sharp e = 1-g.
\]
We obtain the result from knowing that all the quantum integers
are units over $\mathcal{O}_{p}$(\cite[Lem.\,4.1]{GMW04},
Corollary(\ref{c:unit})), and using the definition of $(\ , \
)_{\Sigma}$.

\end{proof}
\begin{defn}
Let $[i]_{p}$ denote
$\frac{A_{p}^{2i}-A_{p}^{-2i}}{A_{p}^{2}-A_{p}^{-2}}$. This is
called the $i$-th quantum integer.
\end{defn}
We can describe the modules $\mathcal{S}_{p}(\Sigma)$ in terms of
`mixed graph' notation in a fixed connected 3-manifold $M$ whose
boundary is $\Sigma$. By a mixed graph, we mean a $p$-admissibly
trivalent graph whose simple closed curves may be colored
$\omega_{p}$ or an integer from the set $\{0,1,\ldots,p-2\}$ where
\[
\omega_{p} = \mathcal{D}_{p}^{-1}\sum_{i=0}^{d_{p}-1}(-1)^{i}[i+1]_{p}  e_{i}.
\]

Using the surgery axiom (S2) in \cite{BHMV3}, we can choose this
fixed 3-manifold to be a handlebody whose boundary is $\Sigma$.
Thus we have
\begin{prop}\label{p:mixed}
A mixed graph in a connected 3-manifold with boundary $\Sigma$
represents an element in $\mathcal{S}_{p}(\Sigma)$. Moreover,
$\mathcal{S}_{p}(\Sigma)$ is generated over $\mathcal{O}_{p}$ by
all the elements given by a mixed graph in a fixed handlebody
whose boundary is $\Sigma$ with the same genus.
\end{prop}
\begin{proof}
The first statement follows from that fact that $V_{p}$ satisfies
the second surgery axiom. The second statement follows from the fact that
every 3-manifold with boundary $\Sigma$ is obtained by a sequence
of 2-surgeries to a handlebody of the same boundary and the
definition of $\mathcal{S}_{p}(\Sigma)$.
\end{proof}

%---------------------------------------------------------------

\section{The first basis for $\mathcal{S}_{p}(S^1\times S^{1})$}

 In this section, we assume $r$ is an odd prime and $p=2r$.
We give a standard basis for $\mathcal{S}_{p}(S^{1}\times S^{1})$.
We need the following lemma before we state our basis.
\begin{defn}
Let $\mu_{i}$ be the eigenvalue for the eigenvector $e_{i}$ of the
twist map on the Kauffman skein module of the solid torus. It is
known in \cite{BHMV1} that $\mu_{i}=(-1)^{i}A^{i(i+2)}$.
\end{defn}
\begin{lem}
For $i\neq j$, we have $\mu_{i}-\mu_{j}$ is equivalent to one of
the following three cases up to a unit in $\mathcal{O}_{p}$.
\begin{enumerate}
    \item 1 if $i\not\equiv j \pmod 2$ and $(j-i)(i+j+2)\not\equiv
    0  \pmod{r}$.
    \item $\sqrt{2}$ if $i\not\equiv j  \pmod{2}$ and $(j-i)(i+j+2)\equiv
    0   \pmod{r}$.
    \item $1+\alpha_{p}^{4}$ if $i\equiv j  \pmod{2}$.
\end{enumerate}
\end{lem}
\begin{proof}
Without loss of generality  we can assume $0\leq i < j \leq
d_{p}-1$. We have

\begin{align*}
\mu_{i}-\mu_{j}& =(-1)^{i}\alpha_{p}^{2i(i+2)}-(-1)^{j}\alpha_{p}^{2j(j+2)}\\
&=(-1)^{i}\alpha_{p}^{2i(i+2)}(1-(-1)^{j-i}\alpha_{p}^{2j^2+4j-2i^2-4i})\\
&\sim 1-(-1)^{j-i}\alpha_{p}^{2(j-i)(i+j+2)}.
\end{align*}
Now we have three cases:

\begin{enumerate}
    \item The hypothesis implies, \[
    \mu_{i}-\mu_{j}\sim  1-(-\alpha_{p}^{2(j-i)(i+j+2)}) \]
    is a unit by Proposition (\ref{p:unit}), as $-\alpha_{p}^{2(j-i)(i+j+2)}$ has
    order divisible by two distinct primes.
    \item The hypothesis implies, \[
    \mu_{i}-\mu_{j}\sim 1\pm i \sim \sqrt{2}.\]
    \item Finally the hypothesis implies that for some $k\leq
    r-1$,
    \[
    \mu_{i}-\mu_{j} \sim 1-\alpha_{p}^{8k} = (1-\alpha_{p}^{8})(1+\alpha_{p}^{8}
    +\dots+\alpha_{p}^{8(k-1)})
    \sim (1-\alpha_{p}^{4})(1+\alpha_{p}^{4})\sim(1+\alpha_{p}^{4}).
    \] As $(1+\alpha_{p}^{8}+\dots+\alpha_{p}^{8(k-1)})$, and
    $(1-\alpha_{p}^{4})$ are units by \cite[Lem 3.1(ii)]{GR97}, and
    Proposition (\ref{p:unit}) respectively.
\end{enumerate}
\end{proof}
\begin{prop}
The number of pairs $(i,j)$ with $0\leq i<j \leq r-2$ such that,
\begin{enumerate}
\item The second case of the above lemma holds is
$(\frac{r-1}{2})$.
\item The third case of the above lemma holds
is $(\frac{r-1}{2})(\frac{r-3}{2})$.
\end{enumerate}
\end{prop}
\begin{proof}
To prove the first part, we look at all pairs $(i,j)$ with
$(j-i)(i+j+2)\equiv 0 \ \pmod r$ which automatically will satisfy
$i\not\equiv j \pmod{2}$. This implies that $i+j+2 = r$. So we
have $(\frac{r-1}{2})$-pairs of such $(i,j)$. Hence the first part
follows. Now for every $0\leq i \leq r-4$, there are
$(\lceil\frac{r-3-i}{2}\rceil)$-$j$'s such that $i\equiv j \pmod
2$. Hence, we have $2(1+2+\ldots+\frac{r-3}{2})$-pairs of such
$(i,j)$. Hence the second part follows.
\end{proof}

\begin{thm}\label{t:first}
Let $\mathcal{B}_{p} = \{ t^{i}(\omega_{p}) \ | \ 0\leq i \leq
d_{p}-1 \}$. Then $\mathcal{B}_{p}$ is a basis for
$\mathcal{S}_{p}(S^{1}\times S^{1})$.
\end{thm}
\begin{proof}
We have
\[ t^{j}(\omega_{p})=\mathcal{D}_{p}^{-1}\sum_{i=0}^{d_{p}-1} (-1)^{i}[i+1]
\mu_{i}^{j}e_{i}.\]

Let $W$ be the matrix which expresses $\mathcal{B}_{p}$ in terms
of $\{e_{0},e_{1},\dots,e_{d_{p}-1}\}$.  The determinant of $W$ is
a unit (product of $(-1)^{i}[i+1]$, see Corollary (\ref{c:unit}))
times $\mathcal{D}_{p}^{-d_{p}}$ times the determinant of the
Vandermonde matrix $[\mu_{i}^{j}]$, where $0\leq i,j\leq
d_{p}-1=r-2$. By the previous lemma
\[
\det
[\mu_{i}^{j}]=\pm\prod_{i<j}(\mu_{i}-\mu_{j})\sim\sqrt{2}^{(\frac{r-1}{2})}
(1+\alpha_{p}^{4})^{(\frac{r-1}{2})(\frac{r-3}{2})} .\] As
$\mathcal{D}_{p}\sim
\sqrt{2}(1+\alpha_{p}^{4})^{(\frac{r-3}{2})}$, we conclude

\begin{equation}\label{e:det}
\det W \sim
(\sqrt{2})^{-(\frac{r-1}{2})}(1+\alpha_{p}^{4})^{-(\frac{r-1}{2})(\frac{r-3}{2})}.
\end{equation}

  Let $\mathcal{W}$ denote the $\mathcal{O}_{p}$-submodule of
$\mathcal{S}_{p}(S^{1}\times S^{1})$ generated by
$\mathcal{B}_{p}$. From the fact that the determinant of $W$ is
non-zero, we conclude that $\mathcal{B}_{p}$ is linearly
independent. Now by Proposition (\ref{p:pairing}), we know
$(e_{i},e_{i})\sim \mathcal{D}_{p}$. Therefore the determinant of
the form (\ref{e:form2}) with respect to this orthogonal basis is
$(\sqrt{2}(1+\alpha_{p}^{4})^{(\frac{r-3}{2})})^{(r-1)}$. By
equation (\ref{e:det}) and the fact $1+\alpha_{p}^{4}\sim
\overline{1+\alpha_{p}^{4}}=1+\alpha_{p}^{-4}$, the determinant of
the form (\ref{e:form2}) with respect to $\mathcal{B}_{p}$ is a
unit. We can conclude that the form on $\mathcal{W}$ is
unimodular. Hence $\mathcal{W} = \mathcal{W}^{\sharp}$ and so the
set $\mathcal{B}_{p}$ forms a basis for
$\mathcal{S}_{p}(S^{1}\times S^{1})$ by equation (\ref{e:dual}).
\end{proof}
\begin{rem}
This theorem and its proof are analogous to
\cite[Thm.\,6.1]{GMW04} and its proof.
\end{rem}

\begin{cor}\label{c:des}

$\mathcal{S}_{p}(\Sigma)$ is generated by 3-manifolds (with no banded
links) with boundary $\Sigma$.
\end{cor}
\begin{proof}
We expand the graph in every element in the Proposition
(\ref{p:mixed}) in terms of linear combinations of banded links
(with some simple curves are colored $\omega_{p}$). Then we
replace any link component (that is not colored $\omega_{p}$) by a
linear combination from the set $\{t^{i}(\omega_{p})| \ 0\leq
i\leq d_{p}-1\}$. Hence the result follows by doing the required
surgery on all the components of the link in every summand.
\end{proof}
\begin{rem}
The above result is true if we replace $p$ by an odd prime as a
corollary of \cite[Thm.\,6.1]{GMW04}.
\end{rem}

%----------------------------------------------------------------

\section{the Frohman Kania-Bartoszynska ideal}
We can apply the results from the previous section to compute the
Frohman Kania-Bartoszynska ideal using the $SU(2)$-theory for
special family of 3-manifolds with torus boundary. Before we do
so, we review this ideal.

\begin{defn}(\cite{FB01})
Let $N$ be a 3-manifold with boundary, we define
$\mathcal{J}_{p}(N)$ to be the ideal generated over
$\mathcal{O}_{p}$ by
\[
\{ I_{p}(M)| \ \text{where}\ M \ \text{is a closed connected
3-manifold containing}\ N \}.
\]
\end{defn}
The importance of this ideal is in being an invariant of
3-manifolds (with boundary) and an obstruction to embedding as
stated in the following propositions.
\begin{prop}(\cite{FB01})
The ideal $\mathcal{J}_{p}$ is an invariant of oriented
3-manifolds with boundary.
\end{prop}
\begin{prop}(\cite{FB01})
If $N_{1}, N_{2}$ are an oriented compact 3-manifolds, and $N_{1}$
embeds in $N_{2}$, then
$\mathcal{J}_{p}(N_{2})\subset\mathcal{J}_{p}(N_{1})$.
\end{prop}
\begin{rem}
Frohman and Kania-Bartoszynska defined this ideal using the
$SU(2)$-TQFT-theory. Afterwards, Gilmer defined this ideal using the
$SO(3)-$TQFT-theory and the 2-theory.
\end{rem}

In general, it is not easy to compute this ideal because we have
infinitely many closed connected 3-manifolds that contains $N$.
Following his work with Masbaum in the case $p$ an odd prime,
Gilmer observed that $\mathcal{J}_{p}(N)$ is finitely generated
based on his result that $\mathcal{S}_{p}(\Sigma)$ is finitely
generated in the case $p$ twice an odd prime as well. We give a
finite set of generators for this ideal for any oriented compact
3-manifold using the $SU(2)-$TQFT-theory which can be obtained by
the following construction.

\begin{defn}\label{d:torus}
Assume $L$ is an ordered link of two components $K$, $J$. Let
$N_{L}$ be the manifold obtained by doing surgery in $S^{3}$ along $K$ in the
complement of $J$.
\end{defn}
\begin{prop}\label{p:generators}
\[
\mathcal{J}_{p}(N_{L})=\langle  I_{p}(M_{i})| \ 0\leq i \leq
d_{p}-1   \rangle,
\] where $M_{i}$ is the 3-manifold obtained by
 doing surgery along the component $K$ and the component
 $J$ with framing $i$ in $S^{3}$.
\end{prop}
\begin{proof}
If $p$ is an odd prime this was proved in \cite{GM04}. With the
help of Theorem (\ref{t:first}), the case $p$ twice an odd prime
follows in the same way.
\end{proof}

%----------------------------------------------------------------

\section{ The Quantization Functors For $p=1,$ and 2}\label{p=2}
In order to understand the relation between $\mathcal{J}_{r}$ and
$\mathcal{J}_{2r}$ when $r$ is an odd prime. We consider the theories
associated to $p=1$ and $p=2$.
\par
We begin by reviewing the quantization functor for $p=1$ in
detail. We start by listing the ring $k_{1} = \mathbb{Z}$, and the
surgery element $\Omega_{1}=\omega_{1}=1$ for this theory defined
in \cite{BHMV1}. We also have $\kappa_{1} = \mathcal{D}_{1} =
\theta_{1}= 1$. One has  $I_{1}(M) = \langle M \rangle_{1} =
\theta_{1}(M) = (-2)^{\sharp k}$ where $\sharp k$ is the number of
components of the banded link in a closed 3-manifold $M$. Then (by
\cite[Prop.\,1.1]{BHMV3}) there exits a unique cobordism generated
quantization functor ($V_{1},Z_{1}$) that extends this invariant.
In fact, this quantization functor can be described explicitly for
surfaces as follows. $V_{1}(\Sigma)$ is the quotient of the
$\mathbb{Z}$-module generated by all 3-manifolds (with banded
links) with boundary $\Sigma$ by the radical of the following form
\[
\langle \ , \ \rangle_{\Sigma}:V_{1}(\Sigma)\times
V_{1}(\Sigma)\rightarrow \mathbb{Z},
\]given by
\[
 \langle[M_{1}],[M_{2}]\rangle_{\Sigma} = \langle
M_{1}\cup_{\Sigma}-M_{2}\rangle_{1} .\] This module is isomorphic
to $\mathbb{Z}$ with any handlebody whose boundary is $\Sigma$ as
a generator. Hence  if $M:\Sigma \rightarrow \Sigma^{'}$,
 then $Z_{1}(M):V_{1}(\Sigma)\rightarrow V_{1}(\Sigma^{'})$ is the just
the multiplication by $(-2)^{\sharp k}$.
\par
 Now we consider the quantization functor for
$p=2$. We start by introducing the ring and its ring of integers
used in this theory
\[ k_{2}=\mathbb{Z}[\alpha_{2},\frac{1}{2}] \quad \text{and}\quad \mathcal{O}_{2} =
\mathbb{Z}[\alpha_{2}].\]

The surgery element for this theory is $\omega_{2}=
\frac{1}{\sqrt{2}}\Omega_{2}$ where $\Omega_{2} = 1 + \frac{z}{2}$
defined in ~\cite{BHMV1}. One has $\mathcal{D}_{2}=\sqrt{2}$, and
$\kappa_{2}=\zeta_{8}$. Therefore the invariant of a closed
connected 3-manifold $M$, which is obtained by doing surgery on
$S^{3}$ along the link $L$, in terms of $\omega_{2}$ is given by
\begin{equation}\label{e:invariant2}
 \langle M\rangle_{2}
= \frac{1}{\sqrt{2}}\kappa_{2}^{\sigma(L)}<L(\omega_{2})>, \
\text{where} < \
> \text{denotes the Kauffman bracket}.
\end{equation}

From this formula, we can easily verify that
\begin{equation}\label{e:sum}
 \langle M_{1}\sharp
M_{2}\rangle_{2} = \sqrt{2}\langle M_{1}\rangle_{2}\langle
M_{2}\rangle_{2}.
 \end{equation}
Now this invariant $\langle M\rangle_{2}$ defined in ~\cite[\S
2]{BHMV3} is involutive and extended to be multiplicative, hence
(by \cite[Prop.\,1.1]{BHMV3}) there exits a unique cobordism
generated quantization functor that extends $\langle M\rangle_{2}$
which is denoted by ($V_{2},Z_{2}$). The modules $V_{2}(\Sigma)$
carry a Hermitian bilinear form defined as follows.
\[
\langle \ , \ \rangle_{\Sigma}:V_{2}(\Sigma)\times
V_{2}(\Sigma)\rightarrow k_{2},
\]given by
\[
 \langle[M_{1}],[M_{2}]\rangle_{\Sigma} = \langle
M_{1}\cup_{\Sigma}-M_{2}\rangle_{2} .\]

 By \cite[1.5 and 6.3]{BHMV3}, $V_{2}(S^{1}\times S^{1})$ is generated
by two elements each of which is a solid torus where the core is
colored either 0 or 1. The pairing in terms of this basis is given
by
\[
\langle 1,1\rangle_{S^{1}\times S^{1}} = \langle S^{1}\times
S^{2}\rangle_{2} = 1.
\]
\[
\langle 1, z\rangle_{S^{1}\times S^{1}} = \langle (S^{1}\times
S^{2},z)\rangle_{2}= \frac{1}{\sqrt{2}}< H > = \frac{1}{2}[2+(-2)]
= 0 = \overline{\langle z,1\rangle}_{S^{1}\times S^{1}}.
\]
Here $H$ is the Hopf link with one of the components is colored
$\omega_{2}$. Finally,
\[
\langle z, z\rangle_{S^{1}\times S^{1}} = \langle (S^{1}\times
S^{2},z\sqcup z)\rangle_{2}= \frac{1}{\sqrt{2}} < K > =
\frac{1}{2}[4+4] = 4.
\]
Here $K$ is the 3-chain link where the middle chain is colored
$\omega_{2}$.

Hence the matrix of the form $\langle \ ,\ \rangle_{S^{1}\times
S^{1}}$ in terms of this basis is given by
$\left(%
\begin{array}{cc}
  1 & 0 \\
  0 &  4 \\
\end{array}%
\right)$.
\par
If we restrict this theory to the category of nonempty connected
objects and connected cobordisms between them, then we have an
integral cobordism theory as before. This follows from the fact
$\left\langle \  \right\rangle_{2}$ is integral as stated in the
proof of \cite[Thm.\,1.1]{GR97} .
\begin{defn}
We define $\mathcal{S}_{2}(\Sigma)$ to be the
$\mathcal{O}_{2}$-submodule of $V_{2}(\Sigma)$ generated by all
connected vacuum states, and we define an
$\mathcal{O}_{2}$-Hermitian bilinear form on
$\mathcal{S}_{2}(\Sigma)$ given by $(\ , \
)_{\Sigma}=\sqrt{2}\langle \ ,\ \rangle_{\Sigma}$.
\end{defn}
\begin{rem}
One could similarly define $\mathcal{S}_{2}^{'}(\Sigma)$ based on
the invariant $\left\langle \ \right\rangle_{2}^{'}$ defined in
\cite[\S.\,1.B]{BHMV3}. In this case, the basis $\{1,z\}$ for
$V_{2}^{'}(S^{1}\times S^{1})$ over $k_{2}$ is also a basis for
$\mathcal{S}_{2}^{'}(S^{1}\times S^{1})$ over $\mathcal{O}_{2}$.
However, this theory is not useful for us in this paper.
\end{rem}

The above basis for $V_{2}(S^{1}\times S^{1})$ does not generate
$\mathcal{S}_{2}(S^{1}\times S^{1})$. The following theorem gives
a basis.

\begin{thm}\label{t:genus1}
Assume that $t$ is the twist map defined in ~\cite{BHMV1}, and
$\mathcal{B} = \{\omega_{2},t(\omega_{2})\}$. Then $\mathcal{B}$
is a basis for $\mathcal{S}_{2}(S^{1}\times S^{1})$, and the form is
unimodular on $\mathcal{S}_{2}(S^{1}\times S^{1})$. Moreover, the matrix of
the form defined in the previous definition in terms of
$\mathcal{B}$ is given by
$\left(%
\begin{array}{cc}
  \sqrt{2} & \frac{1-i}{\sqrt{2}} \\
  \frac{1+i}{\sqrt{2}} & \sqrt{2} \\
\end{array}%
\right)$.
\end{thm}
\begin{proof}
Let $\omega_{2}$ and $t(\omega_{2})$ stands for the elements in
the Kauffman skein module of the solid torus where the core is colored
$\omega_{2}$ and $t(\omega_{2})$ respectively. From the definition
we know that these two elements lie in $\mathcal{S}_{2}(S^{1}\times S^{1})$,
hence $\mathcal{W}=\Span_{\mathcal{O}_{2}}\mathcal{B}\subseteq
\mathcal{S}_{2}(S^{1}\times S^{1})$. The matrix of the form $(\ ,\
)_{S^{1}\times S^{1}}$ is given by
 $\left(%
\begin{array}{cc}
  \sqrt{2} & 0 \\
  0 &  4\sqrt{2} \\
\end{array}%
\right)$, and since the matrix of $\mathcal{B}$ in
terms of $\{1,z\}$ is
given by $\left(%
\begin{array}{cc}
  \frac{1}{\sqrt{2}} & \frac{1}{2\sqrt{2}} \\
  \frac{1}{\sqrt{2}}& \frac{i}{2\sqrt{2}} \\
\end{array}%
\right).$
 Then the matrix $B$ of the form in terms of $\mathcal{B}$ is given by
\[
\left(%
\begin{array}{cc}
  \frac{1}{\sqrt{2}} & \frac{1}{2\sqrt{2}} \\
  \frac{1}{\sqrt{2}}& \frac{i}{2\sqrt{2}} \\
\end{array}%
\right)
\left(%
\begin{array}{cc}
  \sqrt{2} & 0 \\
  0 &  4\sqrt{2} \\
\end{array}%
\right)
\left(%
\begin{array}{cc}
  \frac{1}{\sqrt{2}} & \frac{1}{\sqrt{2}} \\
  \frac{1}{2\sqrt{2}}& \frac{-i}{2\sqrt{2}} \\
\end{array}%
\right)
=\left(%
\begin{array}{cc}
  \sqrt{2} & \frac{1-i}{\sqrt{2}} \\
  \frac{1+i}{\sqrt{2}} & \sqrt{2} \\
\end{array}%
\right) .\] So the form restricted on $\mathcal{W}$ has a unit
determinant. Hence $\mathcal{W} = \mathcal{W}^{\sharp}$. Using
equation (\ref{e:dual}), we get that $\mathcal{W}$ is all of
$\mathcal{S}_{2}(S^{1}\times S^{1})$. In conclusion,
$\{\omega_{2},t(\omega_{2})\}$ is a basis for
$\mathcal{S}_{2}(S^{1}\times S^{1})$.
\end{proof}
\begin{defn}
Let $H_{i_{1}i_{2}\ldots i_{g}}$ be the boundary connected sum of
$g$ solid tori where the core of the $m$-th torus is colored
$i_{m}=0\ \text{or} \ 1$. Also, let
 \[\mathcal{B} =
\{H_{i_{1}i_{2}\ldots i_{g}}| \ (i_{1},i_{2},\ldots,i_{g}) \
\text{is a} \ g\text{-tuple over} \ \{0,1\}\}. \]
\end{defn}
 This set $\mathcal{B}$ is an orthogonal basis for $V_{2}(\Sigma)$, and the pairing
is described as follows:
\begin{prop}
The above set $\mathcal{B}$ forms an orthogonal basis with
respect to the form $\langle \ , \ \rangle_{2}$ given by
\[
\langle H_{i_{1}i_{2}\ldots i_{g}},H_{i_{1}i_{2}\ldots i_{g}}
\rangle_{\Sigma} = 4^{k}(\sqrt{2})^{g-1}, \
\]
where $k = {i_{1}+i_{2}+\ldots +i_{g}}$.
\begin{proof}
By ~\cite[1.5, and 6.3]{BHMV3} $\mathcal{B}$ is a basis. The
result now follows from equation (\ref{e:sum}), and the
computations for $V_{2}(S^{1}\times S^{1})$ after that equation.
\end{proof}
\end{prop}
We can describe $\mathcal{S}_{2}(\Sigma)$ as the
$\mathcal{O}_{2}$-submodule of $V_{2}(\Sigma)$ generated by all
3-manifolds with boundary $\Sigma$ and links sitting inside of
them. As $z=2\sqrt{2}\omega_{2}-2$, one has a similar result to
Corollary (\ref{c:des}) for this theory.
\begin{defn}
Let $H^{'}_{i_{1}i_{2}\ldots i_{g}}$ be the boundary connected sum
of $g$ solid tori where the core of the $m$-th torus is colored
$t^{i_{m}}(\omega_{2})$ for $i_{m}=0$, or 1.

Also, let
\[
\mathcal{B}^{'} = \{H^{'}_{i_{1}i_{2}\ldots i_{g}}| \
(i_{1},i_{2},\ldots,i_{g})\ \text{is a} \ g\text{-tuple over} \
\{0,1\}\}.
\]
\end{defn}
\begin{thm}\label{t:anygenus}
The above set $\mathcal{B}^{'}$ forms a basis for $\mathcal{S}_{2}(\Sigma)$.
\end{thm}
\begin{proof}
Let $(S^{1}\times S^{2})_{ij}$ denote $S^{1}\times S^{2}$ formed
by gluing two solid tori whose cores are colored
$t^{i}(\omega_{2})$, and $t^{j}(\omega_{2})$ where $i,j\in
\{0,1\}$. Let us look at the pairing

\begin{align*}
(H^{'}_{i_{1}i_{2}\ldots i_{g}}, H^{'}_{j_{1}j_{2}\ldots
j_{g}})_{\Sigma} & = \sqrt{2}\langle H^{'}_{i_{1}i_{2}\ldots
i_{g}},H^{'}_{j_{1}j_{2}\ldots j_{g}}\rangle_{\Sigma}\\ & =
\sqrt{2}\langle \sharp_{k=1}^{g}(S^{1}\times S^{2})_{i_{k}j_{k}}\rangle_{2} \\
 &= \sqrt{2}^{g} \prod_{k=1}^{g} \langle
(S^{1}\times S^{2})_{i_{k}j_{k}}\rangle_{2}\\
&=\prod_{k=1}^{g}((S^{1}\times S^{2},t^{i_{k}}(w_{2})\sqcup
t^{j_{k}}(w_{2}))_{2}.
\end{align*}
With a natural order, the matrix of the form in terms of this set
is given by $\bigotimes^{g}B$ ($B$ is defined in the proof of the
previous theorem). This implies that the determinant of this form
is a unit. By a similar argument as in the proof of Theorem
(\ref{t:genus1}), the module generated by this set is all of
$\mathcal{S}_{2}(\Sigma)$.
\end{proof}

We define $I_{2}(M)=\sqrt{2}\langle M_{\flat}\rangle_{2}$ for a
closed 3-manifold $M$ where $\langle M\rangle_{2}$ as defined in
Equation (\ref{e:invariant2}). Also we define the Frohman
Kania-Bartoszynska ideal $\mathcal{J}_{2}$ just as in the previous
section. Now we can compute this ideal easily for all 3-manifolds
using the 2-theory by making use of above results. For example, we
confirm a result of Gilmer and prove it using our basis.
\begin{prop}\cite[Prop.\,15]{G04}\label{p:2ideal}
Let $N_{L}, L$ as defined in Definition (\ref{d:torus}). Also, let
$l$ be the linking number between $K$ and $J$, and $k$ is the
framing of $K$. If $l$ is odd then $\mathcal{J}_{2}(N_{L}) =
\mathcal{O}_{2}$. If $l$ is even, then we have the following:

 \[
 \mathcal{J}_{2}(N_{L}) =
  \begin{cases}
  (\sqrt{2}),
     & \text{if $k \equiv 0  \pmod {4}$};
     \\
     {0}, & \text{if $ k \equiv 2 \pmod {4}$};
    \\
    {\mathcal{O}_{2}}, & \text{if $k \equiv 1 $ or $ 3  \pmod
    {4}$}.
  \end{cases}
  \]
\end{prop}
\begin{proof}
From Theorem (\ref{t:genus1}), we know that
$\mathcal{J}_{2}(N_{L})$ is generated by two elements. In fact, it
is generated by $ I_{2}(M_{i})$ where $M_{i}$ is the 3-manifold
obtained by doing surgery along the component $K$ and the
component $J$ with framing 0 or 1 in $S^{3}$. Now, we use the
formula in \cite[Cor.\,2.4]{BHMV2} to compute the two generators
of this ideal. The computations shows that (note that $\kappa_{2}$
is a unit in $\mathcal{O}_{2}$)
\[
\mathcal{J}_{2}(N_{L}) =
\langle\frac{1}{2}(i^{k+2l}+i^{k}+2),\frac{1}{2}
(i^{k+2l+2}+i+i^{k}+1)\rangle,
\]
Now if we consider all the possibilities, we obtain the required
result.
\end{proof}
Also, we compute this ideal for all 3-manifolds $N_{K}$ that are
obtained by doing surgery on a knot $K$ in the complement of a
tubular neighborhood of an eyeglass graph:$0-0$ in $S^{3}$.

\begin{prop}
Let $l_{1}$ and $l_{2}$ be the linking numbers of $K$ with the
first and the second loops in the eyeglass respectively, and $k$
is the framing of $K$. Then we have
\[
 \mathcal{J}_{2}(N_{K}) =
  \begin{cases}
    {(\sqrt{2})}, & \text{if $l_{1}\equiv l_{2}\equiv 0 \pmod{2}$ and $k\equiv
    0\pmod{4}$};
    \\
    {0}, & \text{if $l_{1}\equiv l_{2}\equiv 0\pmod{2}$ and $ k \equiv 2 \pmod
    {4}$};
    \\
    \mathcal{O}_{2},
     & \text{if any of $l_{1}$, $l_{2}$, $k$ is odd}.
  \end{cases}
\]
\end{prop}
\begin{proof}
Let $m$ be the linking number between the loops. From Theorem
(\ref{t:anygenus}), we know $\mathcal{J}_{2}(N_{K})$ is generated
by four elements. In fact, it is generated by $ I_{2}(M_{i,j})$
where $M_{i,j}$ is the 3-manifold obtained by doing surgery along
the component $K$ and the loops with framing $i, j = 0$ or 1 in
$S^{3}$. As in the proof of the previous proposition, one sees
\begin{multline*}
\mathcal{J}_{2}(N_{K}) =
\langle\frac{1}{2\sqrt{2}}(3+i^{k}+i^{2m}+i^{k+2l_{1}}+i^{k+2l_{2}}+i^{k+2l_{1}+2l_{2}+2m}
),\\ \frac{1}{2\sqrt{2}}
(2+i+i^{k}+i^{2m+1}+i^{k+2l_{1}}+i^{k+2l_{2}+1}+i^{k+2l_{1}+2l_{2}+2m+1}),\\
\frac{1}{2\sqrt{2}}
(2+i+i^{k}+i^{2m+1}+i^{k+2l_{1}+1}+i^{k+2l_{2}}+i^{k+2l_{1}+2l_{2}+2m+1}),\\
\frac{1}{2\sqrt{2}}
(1+2i+i^{k}+i^{2m+2}+i^{k+2l_{1}+1}+i^{k+2l_{2}+1}+i^{k+2l_{1}+2l_{2}+2m+2})\rangle
.
\end{multline*}
If we take all possibilities, we get the required result.
\end{proof}

%--------------------------------------------------------------------

\section{Relating the $r$-th and $2r$-th Theories}

\textit{From now on, we assume that $r$ is an odd prime and
$p=2r$}.
\begin{rem}
The results of this section are slight variations of results of
\cite[\S.\,6]{BHMV3} and \cite[\S.\,2]{BHMV2}. The ring $k_{p}$ is
not exactly the same as the ring denoted this way in \cite{BHMV3}.
\end{rem}
The ring $k_{p}$ will be considered as a $k_{2}$ (or a
$k_{r}$)-module via the homomorphisms defined below. The following
is a slight variation of the maps defined in \cite[\S.\,2]{BHMV2}.
\begin{lem}\label{l:maps}
There are well-defined ring homomorphisms $i_{r}:k_{2}\rightarrow
k_{p}$, $j_{r}: k_{r}\rightarrow k_{p}$  given by

 \[
i_{r}(\alpha_{2}) = \alpha_{p}^{r^{2}}, \ j_{r}(\alpha_{r}) =
 \alpha_{p}^{1+r^{2}} \text{for}\ r\equiv 1\pmod{4}, \text{and}\] \[
 j_{r}(A_{r})=A_{p}^{1+r^{2}} \text{for} \ r\equiv -1\pmod{4}.
 \]

\end{lem}
We need the following remark to prove that these maps are
well-defined.
\begin{rem}
If $\alpha$ is a primitive $n$-th root of unity, then $\alpha^{m}$
is a primitive $\frac{n}{\gcd(n,m)}$-th root of unity.
\end{rem}

\begin{proof}
To prove that the map $i_{r}$ is a well-defined ring homomorphism,
we show $\alpha_{p}^{r^{2}}$ is a primitive 8-th root of unity.
This is true, as $gcd(8r,r^{2})=r$ and $\alpha_{p}$ is a primitive
$8r$-th root of unity. Similarly for $j_{r}$ but we consider two
cases:
\begin{enumerate}
\item For $r\equiv 1 \pmod{4}$, we have $\alpha_{p}^{1+r^{2}}$ is
a primitive $4r$-th root of unity, as $gcd(8r,1+r^{2})=2$ and
$\alpha_{p}$ is a primitive $8r$-th root of unity.

\item For $r\equiv -1 \pmod{4}$, we have $A_{p}^{1+r^{2}}$ is a
primitive $2r$-th root of unity, as $gcd(4r,1+r^{2})=2$ and
$A_{p}$ is a primitive $4r$-th root of unity.

\end{enumerate}

\end{proof}

\begin{cor}\label{c:unit}
The quantum integers $[i]_{p}, 1\leq i\leq d_{p}$ are units in
$\mathcal{O}_{p}$.
\end{cor}
\begin{proof}
We know that the quantum integers $[i]_{r}$  for $1\leq i\leq r-1$ are
units in the $\mathcal{O}_{r}$ see \cite[Lem.\,4.1(iii)]{GMW04}
and \cite[Lem.\,3.1(ii)]{GR97}. So we conclude that $[i]_{p}$ are
units for $1 \leq i \leq r-1 = d_{p}$, as
\[ j_{r}([i]_{r})=(-1)^{i}[i]_{p}, \ \forall \ 1\leq i\leq d_{p},\]
since\[ j_{r}(A_{r}^{2})=-A_{p}^{2}.\]
\end{proof}

 Given any $k_{2}$ (or $k_{r}$)-module, we can define a
$k_{p}$-module by tensoring the original module with $k_{p}$ over
$k_{2}$(or $k_{r}$) respectively. We let $\widehat{V}_{2}(\Sigma)$
(or $\widehat{V}_{r}(\Sigma)$) be the $k_{p}-$module obtained in
this way. We give a relation between
${V}_{1},\widehat{V}_{2},\widehat{V}_{r},$ and ${V}_{2r}$ for any
surface $\Sigma$, but before that we need the following slight
reformulation of \cite[Thm.\,2.1]{BHMV2}.
\begin{thm}\label{t:invariant}
For any closed 3-manifold $M$ with possibly a banded link sitting
inside of it we have,
\begin{equation}\label{e:newformula}
 I_{1}(M) I_{2r}(M) =
i_{r}( I_{2}(M))j_{r}( I_{r}(M)).
\end{equation}
\end{thm}
\begin{proof}
Theorem (2.1) in \cite{BHMV2} states the following:

\begin{equation}\label{e:formula}
\theta_{1}(M)\theta_{2r}(M)
=i_{r}(\theta_{2}(M))j_{r}(\theta_{r}(M)).
\end{equation}

Letting $M = S^{3}$, we get
\[
\mathcal{D}_{2r}^{-1} =
i_{r}(\mathcal{D}_{2}^{-1})j_{r}(\mathcal{D}_{r}^{-1}),
\]

as $\theta_{1}(S^{3})=1$. Now multiply both sides of equation
(\ref{e:formula}) by $\mathcal{D}_{2r}^{-\beta_{1}(M)}$, and
replace it by
$i_{r}(\mathcal{D}_{2}^{-\beta_{1}(M)})j_{r}(\mathcal{D}_{r}^{-\beta_{1}(M)})$
in the right hand side. Then the result follows from the relation
between $\theta$ and $I$.
\end{proof}

We let $\kappa_{n}$ to be an element that plays the role of
$\kappa^{3}$ in \cite{BHMV3}. We define this element as follows:
\[
\kappa_{n}=
\begin{cases}
 \alpha_{4n}^{-6-\frac{n(n+1)}{2}},  & \text{if $n$ is an odd prime};
\\
-\alpha_{4n}^{-6-\frac{n(n+1)}{2}}, & \text{if $n$ is twice an odd
prime}.
\end{cases}
\]
Changing the weight by one multiplies the invariant $\langle \
\rangle_{n}$ by $\kappa_{n}$.
\begin{lem}\label{l:kappa}
For the above ring homomorphisms. We have
\[
\kappa_{p}=i_{r}(\kappa_{2})j_{r}(\kappa_{r})
\]
\end{lem}
\begin{proof}
We have
\[
i_{r}(\kappa_{2})j_{r}(\kappa_{r})=i_{r}(\alpha_{2})j_{r}
(\alpha_{r}^{-6-\frac{r(r+1)}{2}})=\alpha_{p}^{r^{2}}
(\alpha_{p}^{-6-\frac{r(r+1)}{2}})^{1+r^{2}}=-\alpha_{p}^{-6-\frac{2r(2r+1)}{2}}=\kappa_{p}
,\]as
\[
r^{2}+(-6-\frac{r(r+1)}{2})(1+r^{2})\equiv
-6-\frac{2r(2r+1)}{2}+4r \pmod{8r}
\]

\end{proof}

We are now able to give the proof of a result used in $\S2$.
\begin{proof}[Proof of Lemma \ref{l:d}]
The first case follows from \cite[Lem.\,4.1(ii)]{GMW04}. The
second case follows from the following facts from the proof of
Corollary (\ref{c:unit}) and from Theorem (\ref{t:invariant}).
\begin{enumerate}
    \item $\mathcal{D}_{p}=\mathcal{D}_{2r}=i_{r}(\mathcal{D}_{2})j_{r}(\mathcal{D}_{r})$,
    where $\mathcal{D}_{2}=\sqrt{2}$.
    \item $j_{r}(A_{r}^{2})=-A_{2r}^{2}=-\alpha_{2r}^{4}$.
\end{enumerate}
\end{proof}
\begin{thm}\label{t:tensor}
There is a natural $k_{p}$-isomorphism $F:{V}_{1}(\Sigma)\otimes
{V}_{p}(\Sigma)\rightarrow \widehat{V}_{2}(\Sigma)\otimes
\widehat{V}_{r}(\Sigma)$ such that
\begin{equation}\label{e:map}
F([M]_{1}\otimes[M]_{p}) = [M]_{2}\otimes[M]_{r},
\end{equation}
where $M$ is a 3-manifold with banded link (but not linear
combination of links) sitting inside of it.
\end{thm}

\begin{cor}\label{c:induced}
The map in the previous theorem defines a $k_{p}-$isomorphism
between ${V}_{p}(\Sigma)$, and $\widehat{V}_{2}(\Sigma)\otimes
\widehat{V}_{r}(\Sigma)$.

\end{cor}

To prove this theorem, we use the following version of
\cite[Lemm.\,6.4]{BHMV3}.
\begin{lem}
Let $\mathcal{V}$, $W$ be free modules over an integral domain $R$
(with involution) equipped with Hermitian sesquilinear forms
$\langle \ , \ \rangle_\mathcal{V},\ \langle \ , \ \rangle_{W},$
and let $F:\mathcal{V}\rightarrow W$ be a form-preserving linear
map. Let $(V,\langle \ , \ \rangle_{V})$ be the quotient of
$\mathcal{V}$ by the radical of $\langle \ , \
\rangle_{\mathcal{V}}$. Suppose that $\langle \ , \ \rangle_{W}$
is non-degenerate. Suppose $V$ and $W$ are free of finite rank 
and $\langle \ , \ \rangle_{V}$ is
unimodular and furthermore that $\rank(W)\leq \rank(V)$. Then $F$
induces an isometry $\langle \ , \ \rangle_{V}\rightarrow \langle
\ , \ \rangle_{W}$.
\end{lem}
\begin{proof}[ Proof of Theorem \ref{t:tensor}]
It follows from Theorem (\ref{t:invariant}) and Lemma
(\ref{l:kappa}) that formula (\ref{e:map}) defines a
form-preserving linear map. We know already that the form on
$\widehat{V}_{2}(\Sigma)\otimes \widehat{V}_{r}(\Sigma)$ is
non-degenerate. Finally, we have two cases namely,
\begin{itemize}
    \item If $r = 1$, then $F$ is just flipping the tensors. Hence, it
    is an isometry
    \item If $r \geq 3$, then the result follows from the fact
    that $\rank({V}_{1}(\Sigma)\otimes {V}_{2r}(\Sigma)) =
    \rank(\widehat{V}_{2}(\Sigma)\otimes \widehat{V}_{r}(\Sigma))$,
    and the second part of the lemma.
\end{itemize}
\end{proof}

%---------------------------------------------------------------------

\section{The second basis for $ \mathcal{S}_{p}(S^{1}\times S^{1})$}

We give new basis for $\mathcal{S}_{r}(S^{1}\times S^{1})$ that
will be used in constructing another basis for
$\mathcal{S}_{p}(S^{1}\times S^{1})$. To do so, we need the
following lemma.
\begin{lem}
If $0\leq i<j\leq \frac{r-3}{2}$, then the twist coefficients
satisfy $\mu_{i}^{2}-\mu_{j}^{2}\sim \mu_{i}-\mu_{j}\sim
1-A_{r}^{2}$.
\end{lem}
\begin{proof}
Notice that $\mu_{i}=q_{r}^{i^{2}+2i}$ where $q_{r}$ denotes the
primitive $r$-th root of unity given by  $-A_{r}$.
\begin{multline*}
\mu_{i}^{2}-\mu_{j}^{2}=(\mu_{i}+\mu_{j})(\mu_{i}-\mu_{j})=(q_{r}^{i^{2}+2i}+
q_{r}^{j^{2}+2j})(\mu_{i}-\mu_{j})\sim\\
(1+q_{r}^{(j-i)(j+i+2)})(\mu_{i}-\mu_{j})\sim \mu_{i}-\mu_{j}\sim
1-A_{r}^{2}.
\end{multline*}
We used the result of the fourth part of \cite[Lem.\,(4.1)]{GMW04}
in the last equality up to a unit. Also in the one next to last,
we used the fact that $1+q_{r}^{(j-i)(i+j+2)}$ is a unit  by
\cite[Lem.\,(3.1)]{GR97} as $\gcd(r,(j-i)(j+i+2))=1$.
\end{proof}
\begin{thm}\label{t:so3}
Let $\mathcal{B}_{2r} = \{t^{2i}(\omega_{r})|\ 0\leq \ i\leq
d_{r}-1\}$, and $\mathcal{B}_{2r+1}= \{t^{2i+1}(\omega_{r})|\
0\leq \ i\leq d_{r}-1\}$. Then $\mathcal{B}_{2r}$ and $
\mathcal{B}_{2r+1}$ form bases for $\mathcal{S}_{r}(S^{1}\times S^{1})$.
\end{thm}
\begin{proof}
The proof of ~\cite[Thm.\,4.1]{GMW04} now goes through with
$\mu_{i}^{2}$ playing the role of $\mu_{i}$ and $t^{2j}$ playing
the role of $t^{j}$. We use the previous lemma when appropriate to
obtain that $\mathcal{B}_{2r}$ is a basis. To prove that
$\mathcal{B}_{2r+1}$ is a basis, we use the fact that the twist
map $t$ is an isomorphism of the Kauffman skein module of the
solid torus.
\end{proof}
\textbf{Notation:} We use the notation
$\widehat{\mathcal{S}}_{i}(S^{1}\times
S^{1})=\mathcal{S}_{i}(S^{1}\times S^{1})\otimes _{k_{i}} k_{p}$
for $i =2$ or $r$.

\begin{defn}
Let
\[
 \delta_{i} =
  \begin{cases}
  0,
     & \text{if $i + p \equiv \quad 2 \quad or \quad 3  \pmod
     {4}$};
     \\
    1, & \text{if $i + p \equiv \quad 0 \quad or \quad 1  \pmod
    {4}$}.
  \end{cases}
\]
\end{defn}
We defined $\delta_{i}$ so that the following two lemmas
hold.
\begin{lem}\label{l:even}
If $i + \delta_{i} p\equiv 0 \pmod {2}$, then $i + \delta_{i} p
\equiv 0  \pmod {4}$.
\end{lem}
\begin{proof}
We know that $p\equiv 2  \pmod {4}$, now we have two cases to consider
\begin{itemize}
    \item If $\delta_{i} = 0$, then $i + p\equiv 2 \pmod {4}$ as $i$ and $p$ are even. So we
    conclude $i+\delta_{i}p = i \equiv 0  \pmod {4}$.
    \item If $\delta_{i} = 1$, then $i + p\equiv 0  \pmod {4}$ as $i+p$ is even. So we
    conclude $i+\delta_{i}p \equiv i+ p \equiv 0  \pmod {4}$.
\end{itemize}
\end{proof}

\begin{lem}\label{l:odd}
If $i + \delta_{i} p\equiv 1 \pmod {2}$, then $i + \delta_{i} p
\equiv 1 \pmod {4}$.
\end{lem}
A similar proof can be given for this lemma. The following theorem
gives another basis for $\mathcal{S}_{p}(S^{1}\times S^{1})$.

\begin{thm}\label{t:main}
Let $\mathcal{B}_{p} = \{ t^{i+ \delta_{i} p}(\omega_{p}) \ | \
0\leq i \leq d_{p}-1 \}$. Then $\mathcal{B}_{p}$ is a basis for
$\mathcal{S}_{p}(S^{1}\times S^{1})$.
\end{thm}
\begin{proof}

We have that $\Span_{\mathcal{O}_{p}}\mathcal{B}_{p}\subseteq
\mathcal{S}_{p}(S^{1}\times S^{1})$, and  $F(1\otimes
\mathcal{S}_{p}(S^{1}\times S^{1})) \subseteq
\widehat{\mathcal{S}}_{2}(S^{1}\times S^{1}) \otimes
\widehat{\mathcal{S}}_{r}(S^{1}\times S^{1})$ where $F$ is the map
defined in formula (\ref{e:map}). It is enough now to show that
$F(1\otimes \mathcal{B}_{p})$ generates
$\widehat{\mathcal{S}}_{2}(S^{1}\times S^{1}) \otimes
\widehat{\mathcal{S}}_{r}(S^{1}\times S^{1})$ which implies that
$F(1\otimes \mathcal{S}_{p}(S^{1}\times S^{1})) \subseteq$
$\Span_{\mathcal{O}_{p}} F(1\otimes \mathcal{B}_{p})$, i.e
$\mathcal{S}_{p}(S^{1}\times S^{1}) \subseteq
\Span_{\mathcal{O}_{p}} \mathcal{B}_{p}$. Hence, we conclude that
$\mathcal{B}_{p}$ is a basis for $\mathcal{S}_{p}(S^{1}\times
S^{1})$ from the fact that $\rank(\mathcal{S}_{p}(S^{1}\times
S^{1}))= d_{p}$. To prove the claim, let us look at the image of
$\mathcal{B}_{p}$ under $F$.
 \[ \quad t^{i + \delta_{i}
p}(\omega_{p}) \rightarrow t^{i + \delta_{i} p}(\omega_{2})\otimes
t^{i + \delta_{i} p}(\omega_{r}) \quad\] where $i \in
\{0,1,\ldots,d_{p}-1\}$.
\par
Let us consider first all the elements of $\mathcal{B}_{p}$ with
even number of twists, i.e $i+ \delta_{i} p \equiv 0 \pmod {2}$.
By Lemma (\ref{l:even}), we get that $i+ \delta_{i} p \equiv 0
\pmod{ 4}$. Hence those elements get mapped to $t^{4m}(\omega_{2})
= \omega_{2}$ for some $m$, as $t^{4}$ is the identity map in the
2-theory. Also they get mapped to $t^{2j}(\omega_{r})$ for some $0
\leq j \leq d_{p}-1$, as $t^{p}$ is the identity map in the
$SO(3)$-TQFT-theory and $i$ is even. The later elements form the
basis $\mathcal{B}_{2r}$ defined in the previous theorem. In
short, the above elements get mapped to $\omega_{2}\otimes
\mathcal{B}_{2r}$.

Now we consider the elements of $\mathcal{B}_{p}$ with odd number
of twists, i.e $i+\delta_{i}p\equiv 1 \pmod{2}$. By Lemma
(\ref{l:odd}), we get that $i+ \delta_{i} p \equiv 1 \pmod{ 4}$.
Hence those elements get mapped to $t^{4m+1}(\omega_{2})=
t(\omega_{2})$ for some $m$, as $t^{4}$ is the identity map in the
2-theory. Also they get mapped to $t^{2j+1}(\omega_{r})$ for some
$0 \leq j \leq d_{p}-1$, as $t^{p}$ is the identity map in the
$SO(3)$-TQFT-theory and $i$ is odd. The later elements form the
basis $\mathcal{B}_{2r+1}$ defined in the previous theorem. In
short, the above elements get mapped to $\omega_{2}\otimes
\mathcal{B}_{2r+1}$.

Hence the image of $\mathcal{B}_{p}$ under $F$ is a basis for
$\widehat{\mathcal{S}}_{2}(S^{1}\times S^{1})\otimes
\widehat{\mathcal{S}}_{r}(S^{1}\times S^{1})$, i.e generates it as
required.
\end{proof}
\begin{cor}\label{c:isomorphism}
From the above proof, we conclude $\mathcal{S}_{p}(S^{1}\times
S^{1})\cong \widehat{\mathcal{S}}_{2}(S^{1}\times S^{1}) \otimes
\widehat{\mathcal{S}}_{r}(S^{1}\times S^{1})$.
\end{cor}
We do not know if this holds for higher genus surfaces, but it is
clear that $\mathcal{S}_{p}(\Sigma)$ maps into
$\widehat{\mathcal{S}}_{2}(\Sigma) \otimes
\widehat{\mathcal{S}}_{r}(\Sigma)$ under the map $F$.
\begin{prop}\label{p:generators2}
\[
\mathcal{J}_{p}(N_{L})=\langle  I_{p}(M_{i+\delta_{i}p})| \ 0\leq i \leq
d_{p}-1   \rangle,
\] where $M_{i}$ is the 3-manifold obtained by
 doing surgery along the component $K$ and the component
 $J$ with framing $i+\delta_{i}p$ in $S^{3}$.
\end{prop}

Finally, a good question would be: ``Is there a relation between
the Frohman Kania-Bartoszynska ideals in the $SU(2)$- and the
$SO(3)$-TQFT-theories?'' An answer is given by the following
theorem.
\begin{thm}\label{t:ideal}
Let $N$ be an oriented compact 3-manifold with boundary. Then we
have
\[
\mathcal{J}_{p}(N)\subseteq
i_{r}(\mathcal{J}_{2}(N))j_{r}(\mathcal{J}_{r}(N)),
\]
where $i_{r}$ and $j_{r}$ are defined as in the previous section.
Moreover, we have equality if $N$ has a torus boundary.
\end{thm}
\begin{proof}
To prove the inclusion, we have $F(1\otimes
\mathcal{S}_{p}(\Sigma))\subseteq
\widehat{\mathcal{S}}_{2}(\Sigma)\otimes
\widehat{\mathcal{S}}_{r}(\Sigma)$. So if $[M]_{p}\in
\mathcal{S}_{p}(\Sigma)$, then $F(1\otimes [M]_{p}) =
[M]_{2}\otimes [M]_{r}\in \widehat{\mathcal{S}}_{2}(\Sigma)\otimes
\widehat{\mathcal{S}}_{r}(\Sigma)$. By Theorem (\ref{t:invariant})
as $I_{1}(M)=1$, we have
\[ I_{p}(N\cup_{\Sigma}-M)=i_{r}( I_{2} (N\cup_{\Sigma}-M))j_{r}(
I_{r} (N\cup_{\Sigma}-M))\in
i_{p}(\mathcal{J}_{2}(N))j_{p}(\mathcal{J}_{r}(N)).\] So we can
conclude that $\mathcal{J}_{p}(N)\subseteq
i_{p}(\mathcal{J}_{2}(N))j_{p}(\mathcal{J}_{r}(N))$.  Now we prove
the equality in the case of a torus boundary. Let
$M_{i+\delta_{i}p}$ be the solid torus where its core colored
$t^{i+\delta_{i}p}(w)$. From the previous proposition, we have
\begin{align*} \mathcal{J}_{p}(N)&= \langle
I_{p}(N\cup_{S^{1}\times S^{1}}M_{i+\delta_{i}p})
| \ 0\leq i \leq d_{p}-1   \rangle \\
&= \langle i_{r}( I_{2}(N\cup_{S^{1}\times
S^{1}}M_{i+\delta_{i}p})) j_{r} (I_{r}(N\cup_{S^{1}\times
S^{1}}M_{i+\delta_{i}p}))|\ 0\leq i\leq
d_{p}-1\rangle \\
& =i_{r}(\mathcal{J}_{2}(N))j_{r}(\mathcal{J}_{r}(N)).
\end{align*}
The last equality follows from the fact that
$F(1\otimes\mathcal{B}_{p})$ is a basis for
$\widehat{\mathcal{S}}_{2}(S^{1}\times S^{1})\otimes
\widehat{\mathcal{S}}_{r}(S^{1}\times S^{1})$.
\end{proof}

\subsection*{Acknowledgment}
I like to express my deep thanks and gratitude to my Ph.D advisor
Dr. P. Gilmer for all the discussions and the support that were
very important to finish this paper. Also, I like to thank Dr. P.
van Wamelen, and Dr. N. Stoltzfus.

%--------------------------------------------------------------------

\bibliographystyle{amsalpha}

\end{document}